\documentclass{elsarticle}
\usepackage{latexsym,amssymb,amsmath,a4wide}
\usepackage{hyperref}
\newtheorem{theorem}{Theorem}[section]
\newtheorem{lemma}[theorem]{Lemma}

\newtheorem{definition}{Definition}[section]
\newtheorem{assumption}{Assumption}[section]
\newtheorem{remark}{Remark}[section]

\newcommand{\p}{\partial}
\newcommand{\f}{\frac}

\begin{document}

\begin{frontmatter}
\title{\Large Existence of solution to parabolic equations with mixed boundary condition on non-cylindrical domain
%\thanks{}
}
\author[2]{Tujin Kim\corref{1}\fnref{4}}
 \ead{math.tujin@star-co.net.kp}
\author[3]{Daomin Cao \fnref{5}}
 \ead{dmcao@amt.ac.cn}
\cortext[1]{Corresponding author.}
\address[2]{Institute of Mathematics,  Academy of Sciences, Pyongyang, Democratic People's Republic of Korea\fnref{label5}}
\address[3]{Institute of Applied Mathematics,\,AMSS, Chinese Academy of Sciences, Beijing 100190, P. R. China}
\fntext[4]{The author's research was supported by TWAS, UNESCO and
AMSS in Chinese Academy of Sciences.} \fntext[5]{Supported by
Chinese Science Fund for Creative Research Group(10721101)}

\begin{abstract}
  In this paper we are concerned with the initial boundary value problems of linear and semi-linear parabolic equations with mixed boundary
  conditions on non-cylindrical domains in spatial-temporal space. We obtain the existence of a weak solution to the problem. In the case of the linear equation the parts for every type of boundary condition
are any open subsets of the boundary being nonempty the part for
Dirichlet condition at any time. Due to this it is difficult to
reduce the problem to
 one on a cylindrical domain by diffeomorphism of the domain. By a transformation
of unknown function and penalty method we connect the problem to a
monotone operator equation for functions
defined on the non-cylindrical domain. In this way a semilinear problem is considered when the part of boundary for Dirichlet condition is cylindrical.\\
\end{abstract}

\end{frontmatter}

\section{Introduction}

 There are vast literature for parabolic differential
equations on non-cylindrical domain and various methods have been used to study them. In \cite{bhl} the energy inequality for a linear equation with
homogeneous Dirichlet boundary condition is proved, thus unique
existence of solution is studied. For domains expanded along
time existence and uniqueness of solution to initial boundary value
problem of the linear(cf.
 \cite{e}, \cite{ku} and \cite{kp}), semilinear(cf. \cite{krr}) and
nonlinear(cf. \cite{ku}) equations with homogeneous Dirichlet
boundary condition are studied. For such domains and boundary
conditions \cite{krr} also deals with attractor; and \cite{bps}
considers unique existence of solution to a linear
Schr\"odinger-type equation. In \cite{bg} dealing with the Dirichlet
problem, they assume only H\"older continuity on time-regularity of
the boundary. In \cite {ls} semigroup theory is improved and the
obtained result is applied to the initial boundary value problem of
a linear parabolic equation with inhomogeneous Dirichlet condition
on non-cylindrical domain. There are some literatures for unique
existence of initial boundary value problems of linear equations
relying on the method of potentials (see \cite{c} and references
therein). Domains in \cite{l1} and \cite{l2}, where existence,
uniqueness and regularity are studied, are more general, that is,
"initial" condition is given on a hypersurface in spatial-temporal
space instead of the plane $t=0$. In \cite{lms1} optimal regularity
of solution to a special kind of 1-dimensional problem is
considered. Neumann problem of heat equations (cf. \cite{hl}),
parabolic equation with Robin type boundary condition (cf.
\cite{ks}) in non-cylindrical domains  and behavior of solutions to
the initial-boundary value problems of nonlinear equations (cf.
\cite{k} and \cite{t}) are studied. In \cite{lcm} and \cite{mlm}
optimal control and controllability of parabolic equation with
homogenous Dirichlet condition on non cylindrical domain,
respectively are studied.

Also, there are many literatures  for the initial boundary value
problems with mixed boundary conditions. \\Under certain assumptions
the non-cylindrical domains are transformed  to a cylindrical one.
The initial boundary value problems of linear parabolic equations
with mixed time dependent lateral boundary condition on cylindrical
domains are studied (cf. \cite{bfo1}, \cite{bfo2}, \cite{sa}). The
boundary conditions  on the lateral surfaces in \cite{bfo1} may be
either two of the following classical ones: Dirichlet, Neumann and
Robin, but one part for a kind of boundary condition is a connected
and relatively open subset of the lateral surface and the boundary
of the part is tangent to the plane $t=0$. In \cite{bfo2} the first
part is concerned with a classical problem
 on cylindrical domains and the result is applied to the problem with zero initial
condition and lateral mixed boundary conditions on a cylindrical
domain, where the non-cylindrical surface for boundary condition is
transformed to a cylindrical one by a diffeomorphism. The lateral
boundary surface of the cylindrical domain in \cite{bfo2} is also
divided into two parts, and one part $\Gamma_1$ for a kind of
boundary condition is connected and relatively open subset of the
boundary surface and at each point is transverse to the hyperplane
$t=const$. Developing a method in abstract evolution equations,
\cite{sa} is concerned with linear parabolic problems on cylindrical
domains with mixed variable
 inhomogeneous Dirichlet and Neumann conditions. But, here change
along time of distance of the sections of part of boundary for
Dirichlet condition must be dominated by a Lipschitz continuous
function in time $t$.  In \cite{sa} as application of the result,
unique existence of solution to a linear parabolic problem with
homogeneous Dirichlet boundary condition on non-cylindrical domains
is considered.

 Section 4 in
\cite{l3} deals with the linear parabolic problem on non-cylindrical
domain with Robin and Dirichlet boundary conditions, where the
surface for Dirichlet condition is cylindrical type. \cite{lvw} and \cite{lwv}
 study existence, uniqueness and regularity of solutions to the initial
 boundary value problems of linear and semilinear parabolic
 equations on non-cylindrical domains, which is related to the the combustion phenomena. The domains in \cite{lvw} and \cite{lwv} are bordered with a part
 of cylindrical type surface where homogeneous Neumann condition is given, non-cylindrical hypersurfaces where Dirichlet boundary one
is given and planes $t=0,\,\, t=T$. Thus, by change of spatial
independent variable they transform the problems to classical
problems on cylindrical domains where Dirichlet and Neumann
conditions, respectively, are given on cylindrical surfaces. In
\cite{lt} some differential inclusions are studied and the result is
applied to the following problem
$$u_t-\Delta u\in F(u)\,\,\mbox{in}\,\,\Omega,\,\,-\f{\p}{\p n}u\in\beta(u)\,\,\mbox{on}\,\,\gamma,$$
$$u=0\,\,\mbox{on}\,\,\Gamma-\gamma,\,\,u(x,0)=\xi\,\,\mbox{in}\,\,\Omega_0,$$
where $\Omega$ is a non-cylindrical domain in spatial-temporal
space, $\Gamma$ is its lateral surface, $\gamma$ is a part of a
cylindrical surface, $\beta=\p j$ and $j$ is a proper lower
semicontinuous convex function from $R$ to $[0,+\infty]$ with
$j(0)=0$.

On the other hand, in \cite{s} a time-dependent Navier-Stokes
problem on a non-cylindrical domain with a mixed boundary condition
is considered. In \cite{s} the part of boundary where homogeneous
Dirichlet condition is given is cylindrical type and the boundary
condition on the other part of boundary is such a special one that
guarantee existence of solution to an elliptic operator equation
obtained by penalty method.

In this paper we are concerned with linear and semilinear parabolic
equations on non-cylindrical domains with mixed boundary conditions
which may include inhomogeneous Dirichlet, Neumann and Robin
conditions together. In the case of linear equation the parts
 for every type of boundary condition are any open subsets
of the boundary being nonempty the part for Dirichlet condition at
any time. This rises difficulty in reducing the problem to one on
cylindrical domains in \cite{bfo1}, \cite{bfo2} and \cite{sa} or one
in \cite{lvw}, \cite{lwv}, \cite{l3} and \cite{lt}.

 Our idea is to use a transformation of unknown function by
which the problem is connected to a monotone operator equation for
functions defined on the non-cylindrical domain. In this way we can also
consider semilinear equation when the part of boundary for Dirichlet
condition is cylindrical.

This paper is composed of 5 sections. In Section 2 notation, the
problem, the definition of weak solution and the main result are
stated. In Section 3 by a change of unknown function an equivalent
problem is derived. Section 4 is devoted to an auxiliary penalized
problem. In Section 5 the proof of the main result is completed.

\section{Problem and main result}
\setcounter{equation}{0}

Let $\Omega(t)$ be bounded connected domains of $R^N,\,\,\,
Q=\bigcup_{t\in(0,T)}\Omega(t) \times
\{t\},\,\,0<T<\infty,\,\,\,\Sigma=\bigcup_{t\in(0,T)}\p\Omega(t)\times
\{t\}$ and $\Sigma_0,\,\,\Sigma_1$  be open subsets of $\Sigma$ such
that $\bar{\Sigma}_0\cup\bar{\Sigma}_1=\Sigma$ and
$\Sigma_0(t)\equiv\Sigma_0\cap\bar{\Omega}(t)\neq\emptyset\,\forall
t\in(0,T)$. Let $\nu(x,t)$ be outward normal unit vector on the
boundary $\Sigma$ and $n(x,t)$ be outward normal unit
vector on $\p \Omega(t)$ for fixed $t$.\\
 Let $\|y\|_{\Omega(t)}^2=\int\limits_{\Omega(t)}|\nabla y|^2\,dx$ and $|y|_{\Omega(t)}^2=\int\limits_{\Omega(t)}|y|^2\,dx$. Let
$H^1(Q)=W^1_2(Q)$.

 For function $y$ defined on $Q$ define $\beta(y)$
by
$$
\beta(y)=\Big(\int_0^T\|y\|_{\Omega(t)}^2\,dt\Big)^{\f{1}{2}}
$$
whenever the integral make sense. Let
$$
\begin{array}{lc}
D(Q)=\{\varphi:\varphi\in C^\infty(\bar{Q}),\,
\varphi|_{\Sigma_0}=0\},\\
V(Q)=\{\mbox{the completion of $D(Q)$ under the norm}\, \beta(y)\},\\
W(Q)=\{\mbox{the completion of $D(Q)$ in the space}\, H^1(Q)\}
\end{array}
$$
and $\left<\,,\right>$ be duality product between $W(Q)$ and
$W(Q)^*$. By the condition
$\Sigma_0(t)\equiv\Sigma_0\cap\bar{\Omega}(t)\neq\emptyset\,\forall
t\in(0,T)$, $\beta(\cdot)$ is a norm in $D(Q)$.

 We use the following
\begin{assumption}\label{a}
The hypersurface $\Sigma$ belongs to $C^2$  for $x$ and to $C^1$ for
$t$ and for any $t\in[0,T]$ there exist a diffeomorphisms $X(t)$ on
$R^N$ in the class $C^2$ which maps $\Omega(0)$ onto $\Omega(t)$,
$X(0)=I$ and is in $C^1$ for $t$, where $I$ is the unit operator.
\end{assumption}

\begin{remark}\label{r2.1}
Let $\p \Omega(0)\in C^2$ and $\Phi_i(x_1,\cdots,x_N,t)\in
C^{2,1}(R^N\times[0,T]),\,\,i=1,\cdots,N$ be any functions such that
 $\Phi_i(x_1,\cdots,x_N,0)=x_i$ and Jacobian $\f{D\Phi}{Dx}>0$, where $\Phi=\{\Phi_1,\cdots,\Phi_N\}$ and $x=\{x_1,\cdots,x_N\}$.  Then
$\Sigma=\bigcup_{t\in(0,T)}\Sigma(t)\times \{t\}$, where
$\Sigma(t)=\{\Phi(x,t)| x\in\p\Omega(0)\}$ satisfies Assumption 2.1.
\end{remark}

We are concerned with the following initial boundary value problem
\begin{equation}\label{21}
\f{\p y}{\p t} -\sum\limits_{i, j=1}^N \f{\p}{\p
x_i}\Big(a_{ij}(x,t)\f{\p y}{\p x_j}\Big)+\sum\limits_{i=1}^N
b_i(x,t)\f{\p y}{\p x_i}+c(x,t,y)=g(x,t),
\end{equation}

\begin{equation}\label{22}
y|_{\Sigma_0}=\bar{y}|_{\Sigma_0}, \Big(k(x,t)y+\sum\limits_{i,
j=1}^N a_{ij}(x,t)n_i\f{\p y}{\p x_j}\Big)\Big|_{\Sigma_1} =f(x,t),
\end{equation}

\begin{equation}\label{23}
y(x, 0)=y_0(x)\in L_2(\Omega(0)),
\end{equation}
where $a_{ij}(x,t),\,\,b_i(x,t)$ and $c(x,t,r)$ are functions
satisfying the following conditions\\
$(A)\,\,a_{ij}(x,t)\in
W^1_\infty(Q),\,a_{ij}=a_{ji},\,\sum\limits_{i, j=1}^N
a_{ij}(x,t)\xi_i\xi_j\geq\rho|\xi|^2, \, \exists\rho
>0,\,\forall\,\xi \in R^N,\,i=1,\cdots,N,$\\
$(B)\,\, b_i(x,t)\in L_\infty(Q),\,i=1,\cdots,N,$\\
$(C)\,\, c(x,t,r)$ is\,Lipschitz\,continuous\,with\,respect\,to $r$
uniformly for $(x,\,t)$\,and  measurable with respect to $(x,\,t)$
for fixed $r$ and $c(x,t,0)\in L_2(Q)$
and\\
$(D)\,\, \bar{y}\in H^1(Q),\,\,k(x,t)\in
L_\infty(\Sigma_1),\,\,g(x,t) \in L_2(Q),\,\,f(x,t) \in
L_2(\Sigma_1).$
\begin{remark}\label{r2.2}
On a part of $\Sigma_1$ where $k(x,t)=0$ we have Neumann condition.
\end{remark}

 When $y\in C^2(\bar{Q}),\,u \in D(Q)$, in view of
\eqref{22} we have
\begin{equation}\label{27}
\begin{aligned}
\int_Q \f {\p
 y}{\p t}u\,dxdt& =(y(x,T),u(x,T))_{\Omega(T)}-
(y_0,u(x,0))_{\Omega(0)}\\& +\int_{\Sigma_1}yu\cos(\hat{\nu,t})
\,d\sigma-\int_Q y\f {\p u}{\p t}\,dxdt,
\end{aligned}
\end{equation}
\begin{equation}\label{28}
\begin{aligned}
-\int_Q\sum\limits_{i, j=1}^N \f{\p}{\p x_i}\Big(a_{ij}(x,t)\f{\p
y}{\p x_j}\Big)u\, dxdt&=\sum\limits_{i,j=1}^N \int_Q
a_{ij}(x,t)\f{\p y}{\p x_j}\f{\p u}{\p
x_i}dxdt\\&+\int_{\Sigma_1}k(x,t)yu \,d\sigma
-\int_{\Sigma_1}f(x,t)u \,d\sigma,
\end{aligned}
\end{equation}
where $(\hat{\nu,t})$ is the angle between $\nu$ and the positive
direction of $t$-axis. Also, if $y\in L_2(Q)$, then $c(x,t,y)\in
L_2(Q)$ (cf. Lemma 2.2, ch. 2 in \cite{ggz}).

In view of \eqref{27} and \eqref{28}, we introduce
the following
\begin{definition}\label{d2.1}
A function $y$ is called a solution to \eqref{21}-\eqref{23} if $y$
satisfies the following
$$
\begin{aligned}
 &y-\bar{y}\in V(Q),\\
&-\int_Q y\f {\p
 v}{\p t}\,dxdt +\sum\limits_{i,j=1}^N
\int_Q a_{ij}(x,t)\f{\p y}{\p x_j}\f{\p v}{\p
x_i}\,dxdt+\sum\limits_{i=1}^N \int_Q b_i(x,t)\f{\p y}{\p
x_i}v\,dxdt\\
&\hspace{1cm}+\int_Q
c(x,t,y)v\,dxdt+\int_{\Sigma_1}[yv\cos(\hat{\nu,t})+k(x,t)yv]
\,d\sigma\\
&= (y_0,v(x,0))_{\Omega(0)}+\int_Q g(x,t)v\,
dxdt+\int_{\Sigma_1}f(x,t)v\,d\sigma\\
&\hspace{5cm} \forall v \in W(Q)\,\, \mbox{with}\,\, v(x,T)=0.
\end{aligned}
$$
\end{definition}

Our main result of this paper is the following
\begin{theorem}\label{t2.2}
Suppose that conditions $(A),\,(B),\, (C)$ hold and that either
$c(x,t,y)$ is linear with respect to $y$ or $\Sigma_0=\Gamma_0\times
(0,T),\,\,\Gamma_0\subset\p \Omega(0)$ and $\Gamma_0$ is invariant
under the diffeomorphism in Assumption 2.1. Then there exists a
solution to problem \eqref{21}-\eqref{23} provided that $(D)$ is
valid.
\end{theorem}

\section{Transformation of unknown function}
\setcounter{equation}{0}

For the sake of simplicity, we will use the same constants in the
estimates unless confusion will be caused.

It is known that there  exists a function $\psi\in C^2 (\overline{\Omega(0)})$ such
that
$$
\psi(x)>0\,\, \forall x\in \Omega(0),\,\, \psi
|_{\partial\Omega(0)}=0\,\,\mbox{and}\,\, |\nabla\psi|>0 \,\,\forall
x\in \overline{\Omega(0) \backslash \omega_0},
$$
where $\bar{\omega_0}\subset \Omega(0)$ (cf. Lemma 1.1 in
\cite{cik}).
\begin{lemma}\label{l3.1}
There exists a function $\varphi(x,t)\in C^{2,1}(\bar{Q})$ such that
$$\varphi(x,t)>0\,\, \mbox{on}\,\,Q,\,\varphi(x,t)=0\,\,
\mbox{on}\,\,\Sigma\,\,\mbox{and}\,\,-\f {\p \varphi(x,t)}{\p
n}>\eta>0\,\mbox{on}\,\Sigma,
$$
where $C^{2,1}(\bar{Q})$ is the space of functions which are twice
continuously differentiable with respect to $x$ and continuously
differentiable with respect to $t$ on $\bar{Q}$.
\end{lemma}
{\bf Proof} Take $\varphi(x,t)=\psi(X^{-1}(t)x)$, where  $X^{-1}(t)$
is the diffeomorphism from $\Omega(t)$ onto $\Omega(0)$, which is
the inverse of the one given in Assumption \ref{a}. Then the
conclusion follows from the properties of function $\psi$ and
Assumption \ref{a}. $\square$ \vspace*{.3cm}

Let us make a change by
\begin{equation}\label{30}
u=e^{k_1t+k_2\varphi(x,t)}y,
\end{equation}
where $k_1$ and $k_2$ are constants to be determined later. Let
$y\in C^2(\bar{Q})$ and \eqref{22} is satisfied. Then,
$$
\begin{aligned}
&e^{k_1t+k_2\varphi(x,t)}\f{\p y}{\p t}=\f{\p u}{\p t}-(k_1+k_2\f{\p
\varphi}{\p t})u,\\
&-e^{k_1t+k_2\varphi(x,t)}\sum\limits_{i, j=1}^N \f{\p}{\p
x_i}\Big(a_{ij}(x,t)\f{\p y}{\p x_j}\Big)=\\&
\hspace{2cm}-\sum\limits_{i, j=1}^N \f{\p}{\p x_i}\Big(a_{ij}\f{\p
u}{\p x_j}\Big)+2\sum\limits_{ij} a_{ij}k_2\f{\p \varphi}{\p
x_j}\f{\p u}{\p x_i}\\
&\hspace{2cm}+k_2u\Big[\sum\limits_{ij}\f{\p a_{ij}}{\p x_i}\f{\p
\varphi}{\p x_j}-k_2\sum\limits_{ij}a_{ij}\f{\p \varphi}{\p
x_i}\f{\p \varphi}{\p x_j}+\sum\limits_{ij}a_{ij}\f{\p^2 \varphi}{\p
x_i \p x_j}\Big],\\
&e^{k_1t+k_2\varphi(x,t)}\sum\limits_{i=1}^N b_i(x,t)\f{\p y}{\p
x_i}=\sum\limits_{i=1}^N\Big[b_i(x,t)\f{\p u}{\p x_i}-
b_i(x,t)k_2\f{\p \varphi}{\p x_i}u\Big].
\end{aligned}
$$
Also, we have that
$$
\begin{aligned}
 &\sum\limits_{i,j=1}^N
a_{ij}(x,t)n_i\f{\p u}{\p x_j}\Big|_{\Sigma_1}=
\\
&=\Big(\sum\limits_{i,j=1}^N a_{ij}(x,t)n_i\f{\p y}{\p
x_j}e^{k_1t+k_2\varphi(x,t)}+\sum\limits_{i,j=1}^N
a_{ij}(x,t)n_ik_2\f{\p \varphi(x,t)}{\p
x_j}u\Big)\Big|_{\Sigma_1}\\
&=\Big(\sum\limits_{i,j=1}^N a_{ij}(x,t)n_i\f{\p y}{\p
x_j}e^{k_1t+k_2\varphi(x,t)}+k_2\f{\p \varphi(x,t)}{\p
n}\sum\limits_{i,j=1}^N
a_{ij}(x,t)n_in_ju\Big)\Big|_{\Sigma_1}\\
&=\Big(f(x,t)e^{k_1t}-k(x,t)u+k_2\f{\p \varphi(x,t)}{\p
n}\sum\limits_{i,j=1}^N a_{ij}(x,t)n_in_ju\Big)\Big|_{\Sigma_1}
\end{aligned}
$$
where the fact that $\varphi(x,t)=0\,\, \mbox{on}\,\,\Sigma$ and its
corollary $\f{\p \varphi(x,t)}{\p x_j}|_{\Sigma_1}=\f{\p
\varphi(x,t)}{\p n}n_j|_{\Sigma_1} $ have been used.

Taking into account these facts, we have
\begin{equation}\label{31}
\f{\p u}{\p t} -\sum\limits_{i, j=1}^N \f{\p}{\p
x_i}\Big(a_{ij}(x,t)\f{\p u}{\p
x_j}\Big)+\sum\limits_{i=1}^{N}B_i(x,t)\f{\p u}{\p
x_i}+C(x,t,u)=G(x,t),
\end{equation}
\begin{equation}\label{32}
u|_{\Sigma_0}=\bar{u}|_{\Sigma_0},
\Big(K(x,t)u+\sum\limits_{i,j=1}^N a_{ij}(x,t)n_i\f{\p u}{\p
x_j}\Big)\Big|_{\Sigma_1} =F(x,t),
\end{equation}
\begin{equation}\label{33}
u(x, 0)=u_0\equiv y_0(x)e^{k_2\varphi(x,0)}\in L_2(\Omega(0)),
\end{equation}
where
\begin{equation}\label{34}
\begin{aligned}
&B_i(x,t)=b_i(x,t)+2\sum\limits_j a_{ij}k_2\f{\p \varphi}{\p x_j},\\
&C(x,t,u)=e^{k_1t+k_2\varphi(x,t)}c(x,t,
e^{-k_1t-k_2\varphi(x,t)}u)-(k_1+k_2\f{\p \varphi}{\p t})u\\
&\hspace{.6cm}+k_2u\Big[\sum\limits_{ij}\f{\p a_{ij}}{\p x_i}\f{\p
\varphi}{\p x_j}-k_2\sum\limits_{ij}a_{ij}\f{\p \varphi}{\p
x_i}\f{\p \varphi}{\p x_j}+\sum\limits_{ij}a_{ij}\f{\p^2 \varphi}{\p
x_i \p x_j}-\sum\limits_i b_i\f{\p \varphi}{\p x_i}\Big],\\
&G(x,t)=e^{k_1t+k_2\varphi(x,t)}g,\\
&F(x,t)=e^{k_1t}f,\\
&\bar{u}=e^{k_1t+k_2\varphi(x,t)}\bar{y},
\end{aligned}
\end{equation}
and
\begin{equation}\label{35}
K(x,t)=k(x,t)-k_2\f{\p \varphi}{\p
n}\sum\limits_{i,j}a_{ij}n_in_j.
\end{equation}
Now, we take
\begin{equation}\label{35-0}
k_2>0\quad \text{as}\,\, K(x,t)\geq \f{1}{2},
\end{equation}
which is possible by Lemma \ref{l3.1} and $(A)$. Functions
$B_i(x,t),\,C(x,t,u),\,G(x,t),$ $F(x,t),\,\,\bar{u}$ and $K(x,t)$,
respectively, satisfy the conditions for $b_i(x,t),\,c(x,t,r),$
$g(x,t),\,f(x,t),\,\bar{y}$ and $k(x,t)$ in $(B),\,(C)$ and $(D)$.

\begin{lemma}\label{l3.2}
In the sense of Definition \ref{d2.1} existence of solution to
problems \eqref{21}-\eqref{23} and \eqref{31}-\eqref{33} are
equivalent
\end{lemma}
{\bf Proof}  First, let us prove that if $y$ is a solution in the
sense of Definition \ref{d2.1} to problem \eqref{21}-\eqref{23},
then $u$ is a solution in the sense of Definition \ref{d2.1} to
problem \eqref{31}-\eqref{33}.

 For $v\in W(Q)$, put $\bar{v}=e^{-k_1t-k_2\varphi}v$. Then
$$
\begin{aligned}
&-\int_Q y\f {\p
 v}{\p t}\,dxdt=-\int_Q ue^{-k_1t-k_2\varphi}\f {\p
 e^{k_1t+k_2\varphi}\bar{v}}{\p t}\,dxdt\\&
 \hspace*{2.3cm}=-\int_Q u\f {\p
 \bar{v}}{\p t}\,dxdt-\int_Q (k_1+k_2\f{\p \varphi}{\p
 t})uv\,dxdt,\\
&\sum\limits_{i,j=1}^N \int_Q a_{ij}(x,t)\f{\p y}{\p x_j}\f{\p v}{\p
x_i}\,dxdt=\\&\hspace*{1cm} =\sum\limits_{i,j=1}^N \int_Q
a_{ij}(x,t)\f{\p (e^{-k_1t-k_2\varphi}u)}{\p x_j}\f{\p
(e^{k_1t+k_2\varphi}\bar{v})}{\p x_i}\,dxdt\\
&\hspace*{1cm}=\sum\limits_{i,j=1}^N \int_Q a_{ij}\Big[\f{\p u}{\p
x_j}e^{-k_1t-k_2\varphi}-uk_2e^{-k_1-k_2\varphi}\f{\p \varphi}{\p
x_j}\Big]\times\\& \hspace*{4cm}\times\Big[e^{k_1t+k_2\varphi}\f{\p
\bar{v}}{\p x_i}+e^{k_1t+k_2\varphi}k_2\f{\p \varphi}{\p
x_i}\bar{v}\Big]dxdt\\
&\hspace*{1cm}=\sum\limits_{i,j=1}^N \int_Qa_{ij} \Big[\f{\p u}{\p
x_j}\f{\p \bar{v}}{\p x_i}+k_2\f{\p \varphi}{\p x_i}\f{\p u}{\p
x_j}\bar{v}-k_2^2\f{\p \varphi}{\p x_j}\f{\p \varphi}{\p
x_i}u\bar{v}-k_2\f{\p \varphi}{\p x_j}u\f{\p \bar{v}}{\p
x_i}\Big]\,dxdt\equiv I.
\end{aligned}
$$
Integrating by parts in the last term above, we have
$$
\begin{aligned}
I&=\sum\limits_{i,j=1}^N \int_Q a_{ij}(x,t)\f{\p y}{\p x_j}\f{\p
v}{\p x_i}\,dxdt=\sum\limits_{i,j=1}^N \int_Qa_{ij} \f{\p u}{\p
x_j}\f{\p \bar{v}}{\p
x_i}\,dxdt\\&\hspace*{1cm}+\sum\limits_{i,j=1}^N
\int_Q2k_2a_{ij}\f{\p \varphi}{\p x_i}\f{\p u}{\p
x_j}\bar{v}\,dxdt-\int_{\Sigma_1}a_{ij}k_2\f{\p \varphi}{\p
n}n_i n_j u\bar{v}\,d\sigma\\
&\hspace*{1cm}+\sum\limits_{i,j=1}^N \int_Q k_2\Big[\f{\p a_{ij}}{\p
x_i}\f{\p \varphi}{\p x_j}u\bar{v}-k_2a_{ij}\f{\p \varphi}{\p
x_j}\f{\p \varphi}{\p x_i}u\bar{v}+a_{ij}\f{\p^2 \varphi}{\p x_i\p
x_j}u\bar{v}\Big]\,dxdt.
\end{aligned}
$$
Also,
$$
\begin{aligned}
&\sum\limits_i\int_Q b_i\f{\p y}{\p x_i}v\,dxdt=\sum\limits_i\int_Q
b_i\f{\p u}{\p x_i}\bar{v}\,dxdt-\sum\limits_i\int_Q b_ik_2\f{\p
\varphi}{\p x_i}u\bar{v}\,dxdt,\\
&\int_Q c(x,t,y)v\,dxdt=\int_Q
c(x,t,e^{-k_1-k_2\varphi}u)e^{k_1t+k_2\varphi}\bar{v}\,dxdt.
\end{aligned}
$$

 The facts $v\in W(Q)$ and $\bar{v}\equiv e^{k_1t+k_2\varphi}v\in W(Q)$ are
equivalent, and so from above we can see that $u$ is a solution to
\eqref{31}-\eqref{33} in the sense of Definition \ref{d2.1}. In the
same way we can see that if $u$ is a solution to problem
\eqref{31}-\eqref{33} in the sense of Definition \ref{d2.1}, then
$y$ is a solution in the sense of Definition \ref{d2.1} to
\eqref{21}-\eqref{23}. $\square$\vspace{.3cm}

Therefore, in what follows we will consider the existence of a
solution to problem \eqref{31}-\eqref{33}. To this end, in the next
section we will consider an auxiliary problem.

\section{An auxiliary problem}
\setcounter{equation}{0}

The main purpose in this section is to find a function $u^m\in H^1(Q)$ satisfying
the following
\begin{equation}\label{36}
\begin{aligned}
&u^m-\bar{u}\in W(Q),\\
&\int_Q \f{1}{m}\f {\p
 u^m}{\p t}\f {\p
 v}{\p t}\,dxdt-\int_Q u^m\f {\p
 v}{\p t}\,dxdt +
\int_Q \sum\limits_{i,j=1}^Na_{ij}(x,t)\f{\p u^m}{\p x_j}\f{\p v}{\p
x_i}\,dxdt\\&\hspace{1cm} + \int_Q \sum\limits_{i=1}^NB_i(x,t)\f{\p
u^m}{\p x_i}v\,dxdt+\int_Q
C(x,t,u^m)v\,dxdt\\
&\hspace{1cm}+\int_{\Sigma_1}[u^mv\cos(\hat{\nu,t})+K(x,t)u^mv]
\,d\sigma+(u^m(x,T),v(x,T))_{\Omega(T)}\\
&=(u_0,v(x,0))_{\Omega(0)}+\int_Q
G(x,t)v\, dxdt+\int_{\Sigma_1}F(x,t)v\,d\sigma\\
&\hspace*{6cm} \forall v \in W(Q),
\end{aligned}
\end{equation}
where $m$ are positive integers.

We have the following result.
\begin{theorem}\label{t41} Let $k_2$ in \eqref{30} be as \eqref{35-0}.
Then, for some $k_1$ in \eqref{30}, which is taken before
\eqref{312}, there exists a unique solution to problem \eqref{36}.
\end{theorem}
{\bf Proof} Set $u=w+\bar{u}$, define an operator $A_m\in
(W(Q)\mapsto W(Q)^*)$ and an element $L\in W(Q)^*$, respectively, by
$$
\begin{aligned}
&\forall w,v \in W(Q);\\
&\left<A_mw,v\right>=\int_Q \f{1}{m}\f {\p u}{\p t}\f {\p
 v}{\p t}\,dxdt-\int_Q u\f {\p
 v}{\p t}\,dxdt +\sum\limits_{i,j=1}^N
\int_Q a_{ij}(x,t)\f{\p u}{\p x_j}\f{\p v}{\p x_i}\,dxdt\\
&\hspace{1.5cm}+\sum\limits_{i=1}^N \int_Q B_i(x,t)\f{\p u}{\p
x_i}v\,dxdt+\int_Q C(x,t,u)v\,dxdt\\
&\hspace{1.5cm}+\int_{\Sigma_1}\big[uv\cos(\hat{\nu,t})+K(x,t)uv\big]
\,d\sigma+(u(x,T),v(x,T))_{\Omega(T)}
\end{aligned}
$$
and
$$
\left<L,v\right>=(u_0,v(x,0))_{\Omega(0)}+\int_Q G(x,t)v\,
dxdt+\int_{\Sigma_1}F(x,t)v\,d\sigma.
$$

Now, let us consider problem of finding  $w$ such that
\begin{equation}\label{37}
A_mw=L.
\end{equation}
 By the conditions $(A),(B)$ and $(C)$, operator $A_m$ is Lipschitz continuous. For any
$w_1,\,w_2\in W(Q),$ letting $w=w_1-w_2$, we have that
\begin{equation}\label{38}
\begin{aligned}
&\left<A_mw_1-A_mw_2,w\right>=\\&\int_Q \f{1}{m}\f {\p w}{\p t}\f
{\p
 w}{\p t}\,dxdt-\int_Q w\f {\p
 w}{\p t}\,dxdt +\sum\limits_{i,j=1}^N
\int_Q a_{ij}(x,t)\f{\p w}{\p x_j}\f{\p w}{\p x_i}\,dxdt\\
&+\sum\limits_{i=1}^N \int_Q B_i(x,t)\f{\p w}{\p x_i}w\,dxdt +\int_Q
\big[C(x,t,w_1
+\bar{u})-C(x,t,w_2+\bar{u})\big]w\,dxdt\\
&+\int_{\Sigma_1}\big[w^2\cos(\hat{\nu,t})+K(x,t)w^2\big]
\,d\sigma+|w(x,T)|^2_{\Omega(T)}.
\end{aligned}
\end{equation}

On the other hand, by integrating by parts we get
\begin{equation}\label{39}
-\int_Q w\f {\p
 w}{\p
 t}\,dxdt=\f{1}{2}\Big
 [|w(0)|_{\Omega(0)}^2-|w(T)|_{\Omega(T)}^2-\int_{\Sigma_1}w^2\cos(\hat{\nu,t})\,d\sigma\Big].
\end{equation}
From \eqref{38} and \eqref{39} we conclude that for any $w_1,\,w_2\in
W(Q)$
\begin{equation}\label{310}
\begin{aligned}
&\left<A_mw_1-A_mw_2,w\right>=\int_Q \f{1}{m}\f {\p
 w}{\p t}\f {\p
 w}{\p t}\,dxdt+
\int_Q \sum\limits_{i,j=1}^Na_{ij}(x,t)\f{\p w}{\p x_j}\f{\p
w}{\p x_i}\,dxdt\\
&+\int_Q \sum\limits_{i=1}^N B_i(x,t)\f{\p w}{\p x_i}w\,dxdt
 +\int_Q
\big[C(x,t,w_1
+\bar{u})-C(x,t,w_2+\bar{u})\big]w\,dxdt\\
&+\int_{\Sigma_1}\Big[\f{1}{2}w^2\cos(\hat{\nu,t})+K(x,t)w^2\Big]
\,d\sigma+\f{1}{2}|w(0)|_{\Omega(0)}^2+\f{1}{2}|w(x,T)|_{\Omega(T)}^2.
\end{aligned}
\end{equation}
It follows from \eqref{35} and the choice of $k_2$ mentioned above that
\begin{equation}\label{311}
\int_{\Sigma_1}\Big[\f{1}{2}w^2\cos(\hat{\nu,t})+K(x,t)w^2\Big]
\,d\sigma\geq 0.
\end{equation}
Note that $B_i(x,t)$ in \eqref{34} and $K(x,t)$ in \eqref{35} are
independent of $k_1$. Therefore, taking $k_1$ in the expression of $C(x,t,u)$ in
\eqref{34} a negative number small enough independently of $m$, we
have
\begin{equation}\label{312}
\begin{aligned}
&\int_Q \sum\limits_{i,j=1}^N a_{ij}(x,t)\f{\p w}{\p x_j}\f{\p w}{\p
x_i}\,dxdt +\int_Q \sum\limits_{i=1}^N
B_i(x,t)\f{\p w}{\p x_i}w\,dxdt\\
&+\int_Q \big[C(x,t,w_1
+\bar{u})-C(x,t,w_2+\bar{u})\big]w\,dxdt\geq\f{\rho}{2}\int_0^T\|w\|_{\Omega(t)}^2\,dt
\end{aligned}
\end{equation}
By \eqref{310}-\eqref{312} we have
\[
\begin{array}{r}
\left<A_mw_1-A_mw_2,w_1-w_2\right>\geq\alpha\|w\|_{H^1(Q)}^2,\,\,
\exists\,\alpha>0,\,\forall\, w_1,\,w_2\in W(Q)
\end{array}
\]
(Note that $\alpha$ depends on $m$.) Now, by the theory of monotone
operator, there exists a unique solution $w_m$ to problem \eqref{37}
(cf. Theorem 2.2, ch. 3 in \cite{ggz}). Thus, $u^m=w^m+\bar{u}$ is
the solution asserted in the theorem. $\square$

\section{Proof of Theorem \ref{t2.2}}
\setcounter{equation}{0}
 Let $k_1$ be as in the proof of Theorem \ref{t41}, and
$k_2$ as in \eqref{35-0}. When $u^m=w^m+\bar{u}$ is the solution to
\eqref{36} asserted in Theorem \ref{t41}, putting
$w_1=w^m,\,\,w_2=0$, by \eqref{37}, \eqref{38},
\eqref{310}-\eqref{312} we have that
\begin{equation}\label{313}
\begin{aligned}
\int_Q \f{1}{m}\Big|\f {\p
 w^m}{\p t}\Big|^2dxdt+
\int_Q \sum\limits_i\Big|\f{\p w^m}{\p
x_i}\Big|^2dxdt&+|w^m(x,0)|_{\Omega(0)}^2+|w^m(x,T)|_{\Omega(T)}^2\\
&\leq c\left[\left<L,w^m\right>+\left<A_m\bar{u},w^m\right>\right].
\end{aligned}
\end{equation}
Applying Young inequality to the right hand side of \eqref{313} and
taking into account $u^m=w^m+\bar{u}$, we have
\begin{equation}\label{314}
\int_Q \f{1}{m}\Big|\f {\p
 u^m}{\p t}\Big|^2dxdt+
\int_Q \sum\limits_i\Big|\f{\p u^m}{\p
x_i}\Big|^2dxdt+|u^m(x,0)|_{\Omega(0)}^2+|u^m(x,T)|_{\Omega(T)}^2
\leq c,
\end{equation}
where $c$ is independent of $m$.\\

 We claim that for any $v\in W(Q)$
\begin{equation}\label{315}
 \int_Q \f{1}{m}\f {\p u^m}{\p t}\f {\p v}{\p
 t}\,dxdt\rightarrow 0
\,\,\, \mbox{as}\,\,m\rightarrow \infty.
\end{equation}
Indeed, by H\"{o}lder inequality
\[
\Big| \int_Q \f {1}{m}\f {\p u^m}{\p t}\f {\p v}{\p
 t}\,dxdt\Big|\leq\f{1}{\sqrt m}\Big[\int_Q \Big|\f{1}{\sqrt m}\f {\p u^m}{\p
 t}\Big|^2\,dxdt\Big]^{\f{1}{2}}\cdot\Big[\int_Q \Big|\f {\p v}{\p
 t}\Big|^2\,dxdt\Big]^{\f{1}{2}},
\]
which shows \eqref{315} since by \eqref{314} one has that $\int_Q
\f{1}{m}\Big|\f {\p
 u^m}{\p t}\Big|^2dxdt\leq c$.

By \eqref{314} we can also choose a subsequence, which is still
denoted by $\{u_m\}$, such that
\begin{equation}\label{316}
w^m\equiv u^m-\bar{u}\rightharpoonup w,
\mbox\,\,{weakly\,\,in}\,\,V(Q),\,\,\, u_m(T)\rightharpoonup
r\mbox\,\,{weakly\,\,in}\,\,L_2(\Omega(T)).
\end{equation}

First, let $c(x,t,y)$ be linear with respect to $y$, that is, $c(x,t,y)\equiv c(x,t)y$.\\
 Then, $C(x,t,u)$ is also linear with respect to $u$, that is
$C(x,t,u)\equiv C(x,t)u.$ Now put $u\equiv w+\bar{u}$. Then, using
\eqref{315} and \eqref{316} and passing to the limit in \eqref{36},
we have
\begin{equation}\label{317}
\begin{aligned}
&w\equiv\,\,u-\bar{u}\in V(Q),\\
&-\int_Q u\f {\p
 v}{\p t}\,dxdt+
\int_Q \sum\limits_{i,j=1}^Na_{ij}(x,t)\f{\p u}{\p x_j}\f{\p v}{\p
x_i}\,dxdt+\int_Q \sum\limits_{i=1}^NB_i(x,t)\f{\p
u}{\p x_i}v\,dxdt\\
&\hspace{2.4cm}+\int_Q
C(x,t)uv\,dxdt+\int_{\Sigma_1}[uv\cos(\hat{\nu,t})+K(x,t)uv]
\,d\sigma\\
&= (u_0,v(x,0))_{\Omega(0)}+\int_Q G(x,t)v\,
dxdt+\int_{\Sigma_1}F(x,t)v\,d\sigma\\&\hspace{3cm} \forall v \in
D(Q)\,\,\mbox{with}\,\, v(T)=0,
\end{aligned}
\end{equation}
which shows that $u$ is a solution to problem \eqref{31}-\eqref{33}.

Next, let $\Sigma_0=\Gamma_0\times (0,T),\,\,\Gamma_0\subset\p
\Omega(0)$ and $\Gamma_0$ is invariant under the diffeomorphism in
Assumption \ref{a}.\\ Following the method in \cite{s}, we will
prove that the set $\{u^m\}$ of solutions to problem \eqref{36} is
relatively compact in $L_2(Q).$ First, let us prove that the
following two norms in $V(Q)$
\begin{equation}\label{318}
\int_0^T\|w(t)\|_{\Omega(t)}\,dt\,\,
\mbox{and}\,\,\int_0^T\|w(t)\|_{H^1(\Omega(t))}\,dt\,\,\text{are
equivalent},
\end{equation}
where $H^1(\Omega(t))\equiv W^1_2(\Omega(t))$.  It is enough to show
that there exists a constant $C$ independent of $t$ such that
\begin{equation}\label{319}
\|w(t)\|_{H^1(\Omega(t))}\leq C \|w(t)\|_{\Omega(t)}\,\,\, \forall
t\in[0,T]
\end{equation}
Let $w(x)$ be a function defined on $\Omega(t)$ and
$x'\in\Omega(0)$. Then, $w'(x')\equiv w(X(t)x')$ is a function
defined on $\Omega(0)$. By Friedrichs inequality
\[
\int_{\Omega(0)}|w'(x')|^2\,dx'\leq C(\Omega(0),
\Gamma_0)\int_{\Omega(0)}|\nabla w'(x')|^2\,dx'.
\]
Denoting Jacobian of transformation $x'=X^{-1}(t)x$ by $J=\f{D x'}{D
x}$, we have
$$
\begin{aligned}
\int_{\Omega(t)}|w(x)|^2|J|\,dx=\int_{\Omega(0)}|w'(x')|^2\,d
x'&\leq C(\Omega(0), \Gamma_0)\int_{\Omega(0)}|\nabla
w'(x')|^2\,dx'\\
&\leq C(\Omega(0), \Gamma_0)\int_{\Omega(t)}|\nabla
w(x)|^2|J|^{-1}\,dx,
\end{aligned}
$$
where it was considered that in $\nabla w'(x')$ and $\nabla w(x)$  operators $\nabla$ are, respectively, with respect to $x'$ and $x$.

 From this we get
\[
\int_{\Omega(t)}|w(x)|^2\,dx\leq C(t)\int_{\Omega(t)}|\nabla
w(x)|^2\,dx,
\]
where $C(t)$ is continuous in $t\in[0,T]$. This implies \eqref{319}.

On the other hand, by \eqref{314}
\begin{equation}\label{320}
\int_0^T\|w^m\|_{\Omega(t)}^2\,dt\leq c.
\end{equation}
 Let
$\Omega\subset R^N$ such that $\Omega(t)\subset\Omega\,\, \forall
t\in[0,T]$.

For any $m$ let us make $\bar{w}^m(x,t)\in H^1(\Omega\times (0,T))$,
an extension of $w^m(x,t)\in W(Q)$ as follows. Let $w'^m(x')\equiv
w^m(X(t)x')$ on $\Omega(0)$ and denote bounded extensions in
$H^1(R^N)$ by the same(cf. Lemma 1.29, ch. 2 in \cite{ggz}). Then,
$$
\begin{aligned}
\int_{R^N}\big(|w'^m(x')|^2&+|\nabla w'^m(x')|^2\big)\,dx'\leq
c\int_{\Omega(0)}\big(|w'^m(x')|^2+|\nabla w'^m(x')|^2\big)\,dx'\\
&\leq c\int_{\Omega(t)}\big(|w^m(x,t)|^2|J(t)|+|\nabla w^m(x,t)|^2|J(t)|^{-1}\big)\,dx\\
&\leq c\int_{\Omega(t)}\big(|w^m(x,t)|^2+|\nabla
w^m(x,t)|^2\big)\,dx.
\end{aligned}
$$
We take the restriction on $\Omega$ of a function defined by
$\bar{w}^m(x,t)=w'^m(X^{-1}(t)x)$ on $R^N$. Then, we have
$$
\int_\Omega\big(|\bar{w}^m(x,t)|^2+|\nabla
\bar{w}^m(x,t)|^2\big)\,dx \leq c\int_{R^N}\big(|w'^m(x')|^2+|\nabla
w'^m(x')|^2\big)\,dx'.
$$
By \eqref{318}, \eqref{320} and two inequalities above, we get
\begin{equation}\label{321}
\int_0^T\|\bar{w}^m\|_{H^1(\Omega)}^2\,dt\leq c.
\end{equation}
Also, by \eqref{314} we get
\begin{equation}\label{321.1}
\int_{\Omega \times (0,T)} \f{1}{m}\Big|\f {\p
 \bar{w}^m}{\p t}\Big|^2dxdt\leq c,\quad \int_{\Omega \times (0,T)} \f{1}{m}\Big|\f {\p
 \tilde{u}}{\p t}\Big|^2dxdt\leq c,
\end{equation}

\begin{equation}\label{321.2}
 |\bar{w}^m(x,0)|_\Omega\leq c,\quad |\bar{w}^m(x,T)|_\Omega\leq c.
\end{equation}
where $\tilde{u}$ is a bounded extension of $\bar{u}$ and $c$ is
independent of $m$.

 Put $\bar{w}^m(x,t)=0$ for $-T<t<0,\,T<t<2T$. Let
$$
w^m_h(x,t)=\f{1}{h}\int_{t-h}^t\bar{w}^m(x,s)\,ds\quad\mbox{for}\,\,|h|<T.
$$
 Then,
\[
  \f{\p w^m_h(x,t)}{\p
t}=\f{1}{h}\big(\bar{w}^m(x,t)-\bar{w}^m(x,t-h)\big),\,\,\,w^m_h(x,t)|_{\Sigma_0}=0,
\]
which means $w^m_h|_Q\in W(Q).$  Replacing $v$ by
$w^m_h|_Q$ in \eqref{36}, we have
\begin{equation}\label{322}
\begin{aligned}
&\int_Q \f{1}{m}\f {\p
 u^m(x,t)}{\p
 t}\f{1}{h}\big[\bar{w}^m(x,t)-\bar{w}^m(x,t-h)\big]\,dxdt\\
&-\f{1}{h}\int_Q
w^m(x,t)\big[\bar{w}^m(x,t)-\bar{w}^m(x,t-h)\big]\,dxdt\\&-\int_Q
\bar{u}(x,t)\f{\p w_h^m(x,t)}{\p t}\,dxdt + \int_Q
\sum\limits_{i,j=1}^Na_{ij}(x,t)\f{\p u^m}{\p x_j}\f{\p w^m_h}{\p
x_i}\,dxdt\\&+ \int_Q \sum\limits_{i=1}^NB_i(x,t)\f{\p u^m}{\p
x_i}w^m_h\,dxdt+\int_Q C(x,t,u^m)w^m_h\,dxdt\\&
+\int_{\Sigma_1}\big[u^mw^m_h\cos(\hat{\nu,t})+K(x,t)u^mw^m_h\big]
\,d\sigma+(u^m(x,T),w^m_h(x,T))_{\Omega(T)}\\
&=(u_0,w^m_h(x,0))_{\Omega(0)}+\int_Q G(x,t)w^m_h\,
dxdt+\int_{\Sigma_1}F(x,t)w^m_h\,d\sigma.
\end{aligned}
\end{equation}

Assuming $\bar{w}(x,t)\in C^1(\bar{\Omega}\times [0,T])$, let us
estimate
$$I_1\equiv\int\limits_{\Omega\times (0,T)}
\Big|\f{1}{h}\big[\bar{w}(x,t)-\bar{w}(x,t-h)\big]\Big|^2\,dxdt.
$$
First, let $h>0$. Then
\begin{equation}\label{322.1}
\begin{aligned}
I_1&\leq \int_{\Omega\times (h,T)}
\f{1}{h^2}\Big|\int_{t-h}^t\f{\p}{\p
s}\bar{w}(x,s)\,ds\Big|^2dxdt\\
&\hspace*{3cm}+\int_{\Omega\times (0,h)}
\f{1}{h^2}\Big[\bar{w}(x,0)+\int_0^t\f{\p}{\p
s}\bar{w}(x,s)\,ds\Big]^2\,dxdt\\
& \leq \int_{\Omega\times (h,T)}\f{1}{h}\int_{t-h}^t\Big|\f{\p
\bar{w}(x,s)}{\p s}\Big|^2\,ds\,dxdt+
\f{2}{h}|\bar{w}(x,0)|_{\Omega}^2\\
&\hspace*{3cm}+\f{2}{h^2}\int_{\Omega\times
(0,h)}\int_0^t\Big|\f{\p}{\p s}\bar{w}(x,s)\Big|^2\,ds\cdot
h\,dxdt\\
& \leq \f{1}{h}\Big\{(T-h)\int_{\Omega\times (0,T)}\Big|\f{\p }{\p
t}\bar{w}(x,t)\Big|^2\,dxdt+2|\bar{w}(x,0)|_{\Omega}^2\\&
\hspace*{3cm}+2h\int_{\Omega\times (0,T)}\Big|\f{\p}{\p
t}\bar{w}(x,t)\Big|^2\,dxdt\Big\}\\& \leq
\f{1}{h}\Big\{(T+h)\int_{\Omega\times (0,T)}\Big|\f{\p }{\p
t}\bar{w}(x,t)\Big|^2\,dxdt+2|\bar{w}(x,0)|_{\Omega}^2\Big\}.
\end{aligned}
\end{equation}
If $h<0$, using
$$
\begin{aligned}
&I_1\leq \int_{\Omega\times (0,T-|h|)}
\f{1}{h^2}\Big|\int_{t-h}^t\f{\p}{\p
s}w(x,s)\,ds\Big|^2dxdt\\
&\hspace*{1cm}+\int_{\Omega\times (T-|h|,T)}
\f{1}{h^2}\Big[\bar{w}(x,T)+\int_T^t\f{\p}{\p
s}w(x,s)\,ds\Big]^2\,dxdt,
\end{aligned}
$$
 in the same way above we get
\begin{equation}\label{322.2}
\begin{aligned}
&I_1 \leq \f{1}{|h|}\Big\{(T+|h|)\int_{\Omega\times (0,T)}\Big|\f{\p
}{\p t}\bar{w}(x,t)\Big|^2\,dxdt+2|\bar{w}(x,T)|_{\Omega}^2\Big\}.
\end{aligned}
\end{equation}

Since $C^1(\bar{\Omega}\times [0,T])$ is dense in $H^1(\Omega\times
(0,T))$, by \eqref{321.1}, \eqref{321.2}, \eqref{322.1},
\eqref{322.2} for any $\bar{w}^m$ we have
\begin{equation}\label{323}
\begin{aligned}
&\f{1}{\sqrt {m}}\Big\{\int_{\Omega\times (0,T)}
\Big|\f{1}{h}\Big[\bar{w}^m(x,t)-\bar{w}^m(x,t-h)\Big]\Big|^2\,dxdt\Big\}^{\f{1}{2}}\leq\f{c}{\sqrt
{|h|}}.
\end{aligned}
\end{equation}
By \eqref{314} and \eqref{323}, we get
\begin{equation}\label{324}
\begin{aligned}
&\Big|\int_Q \f{1}{m}\f {\p
 u^m(x,t)}{\p
 t}\cdot\f{1}{h}\Big[\bar{w}^m(x,t)-\bar{w}^m(x,t-h)\Big]\,dxdt\Big|\\
& \leq \Big(\int_Q \f{1}{m}\Big|\f {\p
 u^m}{\p
 t}\Big|^2\,dxdt\Big)^{\f{1}{2}}\f{1}{\sqrt {m}}\Big(\int_Q\Big|\f{1}{h}\Big[\bar{w}^m(x,t)-\bar{w}^m(x,t-h)\Big]\Big|^2\,dxdt\Big)^{\f{1}{2}}\\&\leq
\f{c}{\sqrt
{|h|}}. \\
\end{aligned}
\end{equation}

We have
\begin{equation}\label{324.2}
\begin{aligned}
&\Big|\int_Q \sum\limits_{i,j=1}^Na_{ij}(x,t)\f{\p u^m}{\p x_j}\f{\p
w^m_h}{\p x_i}\,dxdt\Big|\leq c\int_0^T\|
u^m\|_{\Omega(t)}\|\f{1}{h}\int_{t-h}^t\bar{w}^m(x,s)\,ds\|_{\Omega(t)}\,dt\\
&\hspace*{2cm}\leq c\int_0^T\| u^m\|_{\Omega(t)}\f{1}{\sqrt
{|h|}}\Big |\int_{t-h}^t\|\bar{w}^m(x,s)\|_{\Omega}^2\,ds\Big
|^{\f {1}{2}}\,dt\leq c/\sqrt {|h|},\\
 &\Big|\int_Q\sum\limits_{i=1}^NB_i(x,t)\f{\p u^m}{\p
x_i}w^m_h\,dxdt\Big|\\
&\hspace*{2cm}\leq c\int_0^T\| u^m\|_{\Omega(t)}\f{1}{\sqrt
{|h|}}\Big|\int_{t-h}^t|\bar{w}^m(x,s)|_{\Omega}^2\,ds\Big|^{\f
{1}{2}}\,dt\leq c/\sqrt {|h|},\\& \Big|\int_Q
C(x,t,u^m)w^m_h\,dxdt\Big|\leq
c\int_Q\Big[|u^m|+|C(x,t,0)|\,\Big]|w^m_h|\,dxdt\leq c/\sqrt {|h|}.
\end{aligned}
\end{equation}
By Assumption \ref{a}, $\left|\f{1}{\sin(\hat{\nu,t})}\right|\geq
\delta>0 $ on $\Sigma_1$. Taking this and the trace theorem into
account, we get
\begin{equation}\label{324.3}
\begin{aligned}&\Big|\int_{\Sigma_1}\big[u^mw^m_h\cos(\hat{\nu,t})+K(x,t)u^mw^m_h\big]
\,d\sigma\Big|\\&\hspace{2cm}=\Big|\int_{\Sigma_1}[u^mw^m_h\cos(\hat{\nu,t})+K(x,t)u^mw^m_h]\f{1}{\sin(\hat{\nu,t})}
\,dxdt\Big|\\&\hspace{2cm}\leq c\int_0^T\|
u^m\|_{\Omega(t)}\f{1}{\sqrt
{|h|}}\Big|\int_{t-h}^t\|\bar{w}^m(x,s)\|_{\Omega}^2\,ds\Big|^{\f
{1}{2}}\,dt\leq c/\sqrt {|h|},
\end{aligned}
\end{equation}
where \eqref{318} and \eqref{320} were used. Also
\begin{equation}\label{325}
\begin{aligned}
&\Big|(u^m(x,T),w^m_h(x,T))_{\Omega(T)}\Big|\leq c\f{1}{\sqrt
{|h|}}\Big|\int_{T-h}^T|\bar{w}^m(x,s)|_{\Omega}^2\,ds\Big|^{\f
{1}{2}}\,dt\leq c/\sqrt {|h|},\\
&\Big|(u_0,w^m_h(x,0))_{\Omega(0)}\Big|\leq c/\sqrt {|h|}.
\end{aligned}
\end{equation}
Similarly to \eqref{324.2}, \eqref{325}, let us estimate $-\int_Q
\bar{u}(x,t)\f{\p w_h^m(x,t)}{\p t}\,dxdt.$
\begin{equation}\label{324.1}
\begin{aligned}
&\Big|-\int_Q
\bar{u}(x,t)\f{\p w_h^m(x,t)}{\p t}\,dxdt\Big|\\
&= \Big|\int_Q \f{\p\bar{u}(x,t)}{\p
t}w_h^m(x,t)\,dxdt-\int_{\Omega(T)}\bar{u}(x,T)w_h^m(x,T)\,dx\Big.\\
&\Big.\hspace{1cm}
+\int_{\Omega(0)}\bar{u}(x,0)w_h^m(x,0)\,dx-\int_{\Sigma_1}\bar{u}(x,t)w_h^m(x,t)\cos(\hat{\nu,t})\,d\sigma\Big|\\
 & \leq c/\sqrt {|h|}.
\end{aligned}
\end{equation}

Also, we get
\begin{equation}\label{326}
\begin{aligned}
&\Big|\int_Q G(x,t)w^m_h\, dxdt\Big|\leq c\int_0^T \f{1}{\sqrt
|h|}\Big|\int_{t-h}^t|\bar{w}^m(x,s)|_{\Omega}^2\,ds\Big|^{\f
{1}{2}}\,dt\leq c/\sqrt {|h|},\\
&\Big|\int_{\Sigma_1}F(x,t)w^m_h\,d\sigma\Big|\leq
c\int_0^T\f{1}{\sqrt
|h|}\Big|\int_{t-h}^t\|\bar{w}^m(x,s)\|_{\Omega}^2\,ds\Big|^{\f
{1}{2}}\,dt\leq c/\sqrt {|h|}.
\end{aligned}
\end{equation}

Let us estimate
$$
I\equiv -\f{1}{h}\int_Q
w^m\big[\bar{w}^m(x,t)-\bar{w}^m(x,t-h)\big]\,dxdt.
$$
Putting $\Omega(t)=\Omega(T)$ for $t>T$,  $\Omega(t)=\Omega(0)$ for
$t<0$ and using $-ab=\f{1}{2}[(a-b)^2-a^2-b^2)]$, we have the
following estimate.

\begin{equation}\label{327}
\begin{aligned}
I&=-\f{1}{h}\int_0^T
\big(\bar{w}^m(x,t),\bar{w}^m(x,t)-\bar{w}^m(x,t-h)\big)_{\Omega(t)}\,dt\\
&=-\f{1}{2h}\int_0^T
\big|w^m(x,t)\big|^2_{\Omega(t)}\,dt+\f{1}{2h}\int_0^T
\big|\bar{w}^m(x,t-h)\big|^2_{\Omega(t)}\,dt\\
&\hspace{3cm}-\f{1}{2h}\int_0^T\big|\bar{w}^m(x,t)-\bar{w}^m(x,t-h)\big|^2_{\Omega(t)}\,dt\\
&=-\f{1}{2h}\int_0^T
\big|w^m(x,t)\big|^2_{\Omega(t)\cap\Omega(t+h)}\,dt-\f{1}{2h}\int_0^T
\big|w^m(x,t)\big|^2_{\Omega(t)\setminus\Omega(t+h)}\,dt\\
&+\f{1}{2h}\int_0^T
\big|\bar{w}^m(x,t-h)\big|^2_{\Omega(t)\cap\Omega(t-h)}\,dt+\f{1}{2h}\int_0^T
\big|\bar{w}^m(x,t-h)\big|^2_{\Omega(t)\setminus\Omega(t-h)}\,dt\\
&\hspace{3cm}-\f{1}{2h}\int_0^T\big|\bar{w}^m(x,t)-\bar{w}^m(x,t-h)\big|^2_{\Omega(t)}\,dt\\
&=-\f{1}{2h}\int_0^T
\big|w^m(x,t)\big|^2_{\Omega(t)\cap\Omega(t+h)}\,dt-\f{1}{2h}\int_0^T
\big|w^m(x,t)\big|^2_{\Omega(t)\setminus\Omega(t+h)}\,dt\\
&+\f{1}{2h}\int_{-h}^{T-h}
\big|\bar{w}^m(x,t)\big|^2_{\Omega(t)\cap\Omega(t+h)}\,dt+\f{1}{2h}\int_0^T
\big|\bar{w}^m(x,t-h)\big|^2_{\Omega(t)\setminus\Omega(t-h)}\,dt\\
&\hspace{3cm}-\f{1}{2h}\int_0^T\big|\bar{w}^m(x,t)-\bar{w}^m(x,t-h)\big|^2_{\Omega(t)}\,dt.\\
\end{aligned}
\end{equation}

Now, put $p=4\,\, \mbox{for}\,\, N=1,\cdots, 4,\,\,
p=\f{2N}{N-2}\,\, \mbox{for}\,\, N\geq5$ and let
$\f{2}{p}+\f{1}{q}=1.$ Applying H\"{o}lder inequality, \eqref{321},
Assumption \ref{a} and the fact that $H^1(\Omega)\hookrightarrow
L_p(\Omega)$, we have
\begin{equation}\label{328}
\begin{aligned}
&\Big|\f{1}{2h}\int_0^T
|\bar{w}^m(x,t-h)|^2_{\Omega(t)\setminus\Omega(t-h)}\,dt\Big|\\
&\leq\f{1}{2|h|}\int_0^T \Big
[\int_{\Omega(t+h)\setminus\Omega(t)}|\bar{w}^m(x,t)|^p\,dx\Big
]^{\f{2}{p}}\cdot\big
[mes\big(\Omega(t+h)\setminus\Omega(t)\big)\,\big]^{\f{1}{q}}\,dt\\
&\leq\f{c}{2|h|}\int_0^T \|\bar{w}^m(x,t)\|^2_{H^1(\Omega)}\,dt\cdot
|h|^{\f{1}{q}}\leq \f{c}{|h|^{(1-1/q)}}.
\end{aligned}
\end{equation}
 Substituting \eqref{328} in the right
hand side of \eqref{327}, we have that
\begin{equation}\label{329}
\begin{aligned}
&-\f{1}{h}\int_Q
w^m(x,t)\big[\bar{w}^m(x,t)-\bar{w}^m(x,t-h)\big]\,dxdt\leq\\
&\hspace{2cm}\f{c}{h^{(1-1/q)}}
-\f{1}{2h}\int_0^T\big|\bar{w}^m(x,t)-\bar{w}^m(x,t-h)\big|^2_{\Omega(t)}\,dt\quad\mbox{if}\,\,h>0,\\
&-\f{1}{h}\int_Q
w^m(x,t)\big[\bar{w}^m(x,t)-\bar{w}^m(x,t-h)\big]\,dxdt\geq\\
&\hspace{1cm}-\f{c}{|h|^{(1-1/q)}}
+\f{1}{2|h|}\int_0^T\big|\bar{w}^m(x,t)-\bar{w}^m(x,t-h)\big|^2_{\Omega(t)}\,dt\quad\mbox{if}\,\,h<0.
\end{aligned}
\end{equation}
Formulas \eqref{320}, \eqref{322},  \eqref{324}-\eqref{326} and
\eqref{329} imply
\begin{equation}\label{330}
\int_0^T\big|\bar{w}^m(x,t+h)-\bar{w}^m(x,t)\big|^2_{\Omega(t)}\,dt\leq
c|h|^{1/q}\quad\mbox{for}\,\, h\in R
\end{equation}
and
\begin{equation}\label{331}
\int_0^h\big|\bar{w}^m(x,t)\big|^2_{\Omega(t)}\,dt+\int_{T-h}^T\big|\bar{w}^m(x,t)\big|^2_{\Omega(t)}\,dt\leq
ch^{1/q}\quad\mbox{for}\,\,h>0.
\end{equation}

 Next, let
\[
\tilde{w}^m(x,t)=\begin{cases} w^m(x,t) &\text{if $(x,\,t)\in Q$}\\
0 &\text{otherwise}.
\end{cases}
\]
Then, when $0<|h|<T$, similarly to \eqref{328} we have
\begin{equation}\label{331.1}
\begin{aligned}
\int_0^T&\big|\bar{w}^m(x,t+h)-\tilde{w}^m(x,t+h)\big|^2_{\Omega(t)}\,dt\\
&\leq
\int_0^T\big|\bar{w}^m(x,t+h)-\tilde{w}^m(x,t+h)\big|^2_{\Omega(t)\setminus\Omega(t+h)}\,dt\\
&\hspace{1cm}+\int_0^T\big|\bar{w}^m(x,t+h)-\tilde{w}^m(x,t+h)\big|^2_{\Omega(t+h)\setminus\Omega(t)}\,dt\\
&\leq
c\int_0^T\big\|\bar{w}^m(x,t+h)\big\|^2_\Omega\,dt\cdot|h|^{\f{1}{q}}\leq
c|h|^{\f{1}{q}}.
\end{aligned}
\end{equation}
From \eqref{330} and \eqref{331.1} we have
\begin{equation}\label{331.2}
\begin{aligned}
&\int_0^T\big|\tilde{w}^m(x,t+h)-w^m(x,t)\big|^2_{\Omega(t)}\,dt
\leq c|h|^{\f{1}{q}}\quad\mbox{for}\,\, 0<|h|<T.
\end{aligned}
\end{equation}

Let $\tilde{h}\in R^N$ and
\[
\Omega_j(t)=\{x\in\Omega(t):\, \mbox{dist}(x,\p
\Omega(t))>2/j\},\quad j=1,2,3, \cdots.
\]
Then, by \eqref{318} and \eqref{320}
\begin{equation}\label{332}
\begin{aligned}
\int_{\Omega(t)\setminus\Omega_j(t)}&|w^m(x,t)|^2\,dx=\\&=\Big[\int_{\Omega(t)\setminus\Omega_j(t)}|w^m(x,t)|^p\,dx\Big
]^{\f{2}{p}}\cdot\big
[mes\big(\Omega(t)\setminus\Omega_j(t)\big)\,\big]^{\f{1}{q}}\,dt\\
&\leq c\|w^m(x,t)\|^2_{\Omega(t)}\cdot\Big( 1/j \Big )^{\f{1}{q}},
\end{aligned}
\end{equation}
where $c$ depends only on $Q$.

 Now, if $|\tilde{h}|<1/j$, then
$x+s\tilde{h}\in \Omega_{2j}(t)$ provided $x\in\Omega_j(t)$ and
$s\in[0,1]$. For $w^m\in C^\infty(\Omega(t))$,
\begin{equation}\label{333}
\begin{aligned}
\int_{\Omega_j(t)}|w^m(x+\tilde{h},t)-w^m(x,t)|^2&\,dx\leq
\int_{\Omega_j(t)}\,dx\Big
[\int_0^1\big|\f{d}{ds}w^m(x+s\tilde{h},t)\big|\,ds\Big]^2\\
&\leq \int_{\Omega_j(t)}\,dx\Big [\int_0^1\big|\nabla
w^m(x+s\tilde{h},t)\big|\cdot\big|\tilde{h}\big|\,ds\Big]^2\\&\leq
|\tilde{h}|^2\int_0^1\int_{\Omega_{j}(t)}\big|\nabla w^m(x+s\tilde{h},t)\big|^2\,dxds\\
&\leq |\tilde{h}|^2\int_{\Omega_{2j}(t)}|\nabla w^m|^2\,dx
\leq(1/j)^2\|w^m\|_{\Omega(t)}^2.
\end{aligned}
\end{equation}
Since $C^\infty(\Omega(t))$ is dense in $H^1(\Omega(t))$,
\eqref{333} is valid for any $w^m\in H^1(\Omega(t))$. By
\eqref{332}, \eqref{333} and \eqref{320} we have that
\begin{equation}\label{334}
\int_Q\big|\tilde{w}^m(x+\tilde{h},t)-\tilde{w}^m(x,t)\big|^2\,dxdt\leq
c(1/j)^{1/q}\quad\mbox{for}\,\,|\tilde{h}|<1/j
\end{equation}
where $c$ is independent of $m$.

From \eqref{331.2} and \eqref{334}
we get
\begin{equation}\label{335}
\int_Q\big|\tilde{w}^m(x+\tilde{h},t+h)-\tilde{w}^m(x,t)\big|^2\,dxdt\rightarrow
0\quad \mbox{as}\,\,(\tilde{h},h)\rightarrow
0\,\,\mbox{in}\,\,R^{N+1}.
\end{equation}
From \eqref{331} and \eqref{332}, we get
\begin{equation}\label{336}
\forall\varepsilon,\, \exists Q_\varepsilon\,\,
\mbox{such\,that}\,\, \bar{Q}_\varepsilon\subset Q:
\int_{Q\backslash
Q_\varepsilon}\big|w^m(x,t)\big|^2\,dxdt<\varepsilon.
\end{equation}
By \eqref{335} and \eqref{336} we know that the set $\{w^m\}$ is
relatively compact in $L_2(Q)$ (cf. Theorem 2.32 in \cite{af}).
Thus, we can choose a subsequence, which is still denoted by
$\{w^m\}$, such that $ w^m\rightarrow w \in L_2(Q)$. Therefore
$C(x,t,u^m)\rightarrow C(x,t,u)$ in $L_2(\Omega(t))$ for a.e. $t$,
where $u\equiv w+\bar{u}$.

Therefore, using \eqref{316} and passing to the limit in \eqref{36},
we have \eqref{317} which shows that $u$ is a solution to problem
\eqref{31}-\eqref{33}. $\square$

 \end{document}